\input amstex
\documentstyle{amsppt}
\magnification 1200
\vcorrection{-1cm}
\NoBlackBoxes
\input epsf

\let\em=\it

\rightheadtext{ Corrigendum to ``A flexible affine $M$-sextic'' }

%

\topmatter
\title
          {\rm Corrigendum to the paper} ``A flexible affine $M$-sextic
          which is algebraically unrealizable''
\endtitle
\author
          S.~Fiedler--Le Touz\'e,
          S.~Orevkov, and E.~Shustin
\endauthor
\abstract
We prove the algebraic unrealizability of certain isotopy type of plane affine
real algebraic $M$-sextic which is pseudoholomorphically realizable.
This result completes the classification up to isotopy of real algebraic affine $M$-sextics.
The proof of this result given in a previous paper by the first two authors was incorrect.
\endabstract
\address
          IMT, l'universit\'e Paul Sabatier, 118 route de Narbonne, 31062 Toulouse, France
\endaddress
\email    severine.fiedler\@live.fr
\endemail
\address
          IMT, l'universit\'e Paul Sabatier, 118 route de Narbonne, 31062 Toulouse, France
and
          Steklov Math. Institut, Gubkina 8, 119991 Moscow, Russia
\endaddress
\email    orevkov\@math.ups-tlse.fr
\endemail
\address  School of Math. Sciences, Tel Aviv Univ., Ramat Aviv, 69978
          Tel Aviv, Israel
\endaddress
\email    shustin\@math.tau.ac.il
\endemail
\thanks
          The second author was partially supported by RFBR grant 17-01-00592a
\smallskip
          The third author has been supported by the
          Israeli Science Foundation grant no. 176/15.
\endthanks
\endtopmatter

\def\figTh  {1}
\def\figone {1}
\def\figTwoABC {2}
\def\figtwod   {3}
\def\figA  {4}
\def\figB  {5}
\def\figError  {6}
\def\figPencil {7}
\def\figSev    {8}

\def\figtwoa {\figTwoABC(a)}
\def\figtwob {\figTwoABC(b)}
\def\figtwoc {\figTwoABC(c)}

\def\sectHRG    {1}
\def\sectResolv {2}
\def\sectError  {3}

\def\thone      {1}
\def\lemone     {1}
\def\lemtwo     {2}
\def\remone     {1}

\def\eqone     {3}
\def\eqtwo     {2}
\def\eqthree   {1}

\def\refB      {1}
\def\refSevOne {2}
\def\refSevTwo {3}
\def\refFO     {4}
\def\refGdan   {5}
\def\refG      {5b}
\def\refIKS    {6}
\def\refKS     {7}
\def\refTopol  {8}
\def\refAFST   {9}
\def\refAa     {10}
\def\refCrelle {11}
\def\refMMJ    {12}
\def\refAaX    {13}
\def\refSh     {14}
\def\refV      {15}

\let\refone   =\refAaX 
\let\refthree =\refMMJ 
\let\reffour  =\refSh  
\let\reffive  =\refIKS
\let\refsix   =\refAaX
\let\refseven =\refG

\def\R{\Bbb R}

\def\RP{\Bbb{RP}}
\def\CP{\Bbb{CP}}
\def\hirz{\Cal F}
\def\eps{\varepsilon}

\def\On{O^{\text{ne}}}
\def\Oe{O^{\text{e} }}

\document
The main theorem of the paper [\refFO] states that
the arrangement $B_2(1,4,5)$ in $\RP^2$ (see Figure \figTh) is
unrealizable by a union of a line and a real smooth algebraic sextic curve.
The precise statement is:

\proclaim{ Theorem 1 } Let $C$ be a real algebraic curve of degree $6$
in $\RP^2$ and $L$ a line. Then there does not exist an ambient isotopy
of $\RP^2$ which deforms $C$ and $L$ into the curve and the line in Figure \figTh.
\endproclaim

\midinsert
\epsfxsize=40mm\centerline{\epsfbox{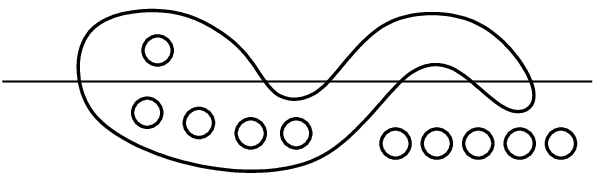}}
\botcaption{ Figure \figTh }
The arrangement $B_2(1,4,5)$
\endcaption
\endinsert

Recently, the first author found a mistake in the
final part of the proof of Theorem 1 given in [\refFO]
(we discuss this mistake in detail in Section \sectError\ below).
However the result is correct and here we give another proof.

An affine smooth irreducible real algebraic curve $A$ in $\R^2$ of degree $d$ is an
{\it affine $M$-curve} if it has maximal possible number of connected components,
which is equal to $g+d$ where $g=(d-1)(d-2)/2$
is the genus of the complexification of $A$.
This condition is equivalent to the fact that the projective closure
of $A$ is an $M$-curve (i.e. it has $g+1$ connected components) and all intersections
with the infinite line are real and transverse and sit on the same connected component of
the closure of $A$. Thus, if Figure \figTh\ were algebraically realizable,
it would provide an affine $M$-sextic in the affine plane $\RP^2\setminus L$.

A classification of affine $M$-sextics up to isotopy was started in [\refKS, \refTopol]
and completed in [\refMMJ, Theorem 1.1] assuming that [\refFO] is correct.
So, here we fill a gap in the proof of this classification as well.
Note that a pseudo-holomorphic classification of affine  $M$-sextics was previously
obtained in [\refTopol], and it differs from the algebraic one. Three arrangements
are realizable pseudo-holomorphically, but not algebraically; see [\refMMJ].
The arrangement $B_2(1,4,5)$ in Figure \figTh\ is one of them.
This is why it is more difficult to exclude it.


We prove Theorem \thone\ arguing by contradiction and proceed in three steps:
\roster
\item"(i)" assuming that a smooth sextic curve $C_0$ arranged
           with respect to the line $L$ as shown in Figure \figone\ exists,
           we derive that there exists a real elliptic sextic curve $C_9$ with
           $9$ nodes located with respect to $L$ as shown in Figure \figtwoa\
           (see Lemma \lemone\ in Section \sectHRG);
\item"(ii)" from the existence of a sextic $C_9$ we derive the existence of an
            elliptic real sextic having $7$ nodes (five isolated and two non-isolated)
            and a singularity $A_3$, and located with respect to $L$ as shown
            in Figure~\figtwob (see Lemma \lemtwo\ in Section \sectHRG);
\item"(iii)" we prohibit the existence of the latter real elliptic sextic using
             a suitable version of cubic resolvent (see Section \sectResolv).
\endroster

\midinsert
\centerline{\epsfxsize=120mm\epsfbox{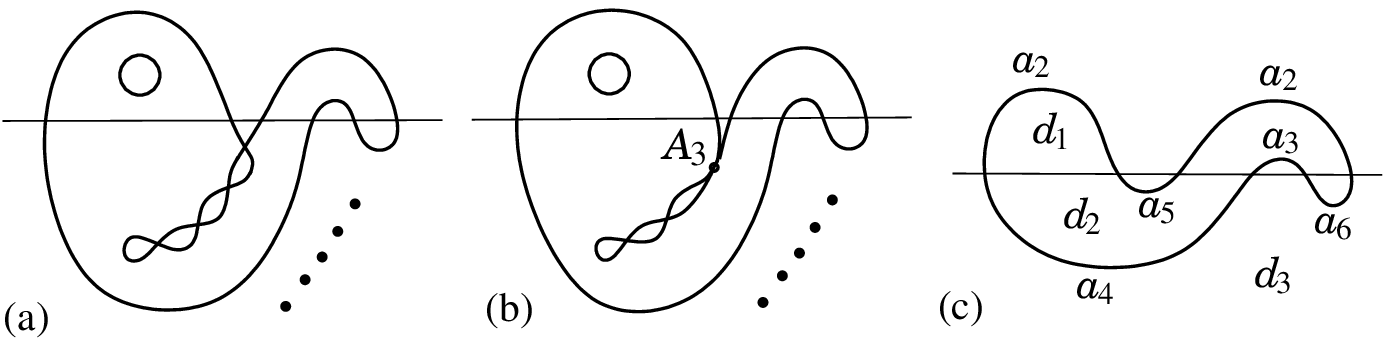}}
\botcaption{ Figure \figTwoABC }
(a) See Lemma \lemone. (b) See Lemma \lemtwo. (c) Notation $a_i$, $d_i$.
\endcaption
\endinsert

So, the general scheme of the proof of Theorem \thone\ is almost the same as for
the proof in [\refMMJ] of algebraic unrealizability of the affine sextic $C_2(1,3,6)$.
However, there is a difference in the last step. In [\refMMJ], the cubic resolvent is
algebraically realizable but its mutual position with respect to the axis is
pseudo-holomorphically unrealizable. Here the situation is opposite: the
mutual position of the resolvent and the axis is pseudo-holomorphically realizable,
but the resolvent itself is algebraically unrealizable.

Note also that only the $A_3$ singularity is needed in Step (iii), thus
we may undo all the remaining $A_1$ singularities obtained in Step (ii).
However, we do not know how to attain the $A_3$ singularity (keeping the
required position of the curve with respect to the line $L$) without passing
through a genus 1 nodal curve.


\head    \sectHRG. Application of Hilbert-Rohn-Gudkov method
\endhead

The following notation will be used in the proof.
Denote the nonempty oval of $C_0$ by $\On$. In what follows, we deform $C_0$
in certain families, and the corresponding non-empty oval will be denoted by $\On$ as well.
If $\On$ degenerates into a loop with a singular point, we continue to use the
same notation for this loop. The open disk bounded by $\On$ will be denoted by $d_0$.
The line $L$ cuts $\On$ into six arcs, which we denote $a_1,\dots,a_6$
according to Figure~\figtwoc.
We also use the notation $d_1,d_2,d_3$ for the three of the connected
components of $\RP^2\setminus(\On\cup L)$ designated in Figure \figtwoc.
By $\Oe$ we denote the empty oval in the domain $d_1$ (it will remain oval in
all further deformations).

\proclaim{Lemma \lemone}
Suppose that there exists a sextic curve $C_0$ shown in Figure~\figone.
Then there exists a real irreducible sextic curve $C_9$ with $9$ nodes
located with respect to the line $L$ as shown in Figure~\figtwoa.
\endproclaim

\demo{ Proof }
We construct the curve $C_9$ inductively. Abusing notation we denote
a plane curve and its defining polynomial by the same symbol.

Start with a pencil of sextic curves $\{C_0^{(t)}=C_0+\eps tK_0^2, \ \eps=\pm1,\ t\ge0\}$,
where $K_0$ is a generic real cubic curve passing through a point $p$
chosen on the arc $a_6$. Choose $\eps$ so that the disk $d_0$ contracts as $t$ grows.
Furthermore, the oval $\On$ always intersects $L$ as shown in Figure~\figone,
the disks bounded by the empty ovals inside $d_0$ grow, while the disks bounded by the
empty ovals inside $d_3$ shrink.
Note that $t$ cannot tend to $\infty$ without degeneration of $C_0^{(t)}$.
Indeed, otherwise the curve $C_0^{(t)}$ would approach the double
cubic curve $K_0^2$.
This, however, is impossible: we consider the real line through
a fixed point inside the oval $\Oe$ and a point embraced by one of the empty
ovals in the domain $d_2$, and then, on any sufficiently close real line,
we will observe a pair of real intersection points
with $C_0^{(t)}$ that approach the intersection point with $L$ (see Figure \figtwod, left).
This finally would yield that $K_0$ contains a segment of $L$. A contradiction.

Due to the general choice of $K_0$, the first degeneration
is a sextic curve $C_1$ with one node (see more detailed arguments in
the proof of the induction step below).
\if01{
Since the family of sextics having a singularity more complicated than a node has
codimension $\ge2$ in the space of sextics passing through $p$, and
Due to the general choice of $K_0$,
the pencil $\{C_0^{(t)}\}_{t\in\R}$ does not hit the families of sextics passing
through $p$ and having singularity collection different from one node, since such families
have codimension $\ge2$ in the linear system of sextics passing through $p$. Hence,
the first degeneration
is a sextic curve $C_1$ with one node.
}\fi 
Note that this node cannot
join $\On$ with the empty oval in the domain $d_1$, since they
form a positive complex oriented injective pair
(see [\refFO, Figure 9]).
\footnote{We refer to [\refV, \S2.1 and \S2.4B] for a definition of the complex
orientations and of positive/negative injective pairs of ovals respectively.
Also [\refV] can be used as a general introduction to the subject.}

Further, the node cannot
join the arc $a_4$ with an empty oval in the domain $d_2$.
Indeed, suppose there exists such a nodal degeneration $C^*$.
Let us consider the pencil of lines through a point inside $\Oe$.
By [\refFO, Lemma 2.2], the lines of this pencil passing through the exterior ovals
do not separate the ovals in $d_2$ from each other and from the arc $a_5$.
Hence (cf. [\refFO, \S6.2] and [\refAa, \S4.5]) the braid of $C^*$ with respect to this pencil coincides
with that of a curve $C^{**}$ obtained from Figure \figtwoa\ by joining one of the ovals
in $d_2$ with $a_4$ through one more node. Thus there exists such a pseudo-holomorphic
curve $C^{**}$ which contradicts the genus formula.
Hence,
\roster
\item"(1)"
           either the node joins a pair of empty ovals in the domain $d_2$,
\item"(2)"
           or the node joins an empty oval with the arc $a_5$,
\item"(3)"
           or the node is an isolated point obtained by shrinking one of
           the empty ovals in the domain $d_3$.
\endroster

For the induction step, suppose that $C_k$, $1\le k\le 8$, is a real irreducible
sextic curve with $k$ real nodes and such that
\roster
\item"(i)" there exists a smoothing of $C_k$ into a smooth sextic
           as shown in Figure \figone,
\item"(ii)" all nodes are real, each non-isolated node joins either a pair of empty
            ovals in the domain $d_2$, or an empty oval in the domain $d_2$ with
            the arc $a_5$, while each isolated node is obtained by shrinking an
            oval in the domain $d_3$.
\endroster

By [\refthree, Proposition 2.4(b)], we can suppose that all $k$ nodes and the point
$p$ are in general position. Consider a pencil
$$
    \{C_k^{(t)}=C_k+\eps tK_k^2,\quad \eps=\pm1,\quad t\ge0\}\ ,
                                                                  \eqno(\eqthree)
$$
where $K_k$ is a generic real cubic passing through the nodes of $C_k$ and the point
$p\in a_6$, and $\eps$ is chosen so that $d_0$ shrinks as $t$ grows.
The argument, used for the base of induction, ensures that there must be a
degeneration at some $t\in(0,\infty)$. We claim that the first degeneration
$C^*$ is a $(k+1)$-nodal curve $C_{k+1}$ possessing the
above properties (i) and (ii). We explain this just in the most difficult case of $k=8$.
By [\refGdan, Theorem 1] we can suppose that all node of $C_8$ are in general position.
Blowing up the $8$ fixed nodes, we obtain curves in the $3$-dimensional linear system
$|D|=|6L-2E_1-...-2E_8|$ on a general del Pezzo surface $\Sigma$ of degree $1$.
By [\reffive, Lemma 9(1)],
the curves that are not immersed form a set
of dimension at most $1$ in $|D|$. Fixing the point $p\in \Sigma$,
we obtain that the pencils spanned by non-immersed curves and the (unique) double
curve passing through $p$, sweep a subset of dimension $\le 2$ in $|D|$,
while the (smooth) blown up curve $C_8$ can be moved to a
general position in $|D|$.
Hence, the considered pencil of sextics (\eqthree)
does not contain non-immersed curves (except for the double cubic).
We then see that the degenerate curve $C^*$ cannot have an
immersed singularity more complicated than a node by the genus formula and the Harnack-Klein bound
stating that the number of connected components of the real point set of the normalization of the curve does not exceed genus plus one.
For the same reason, an extra singularity cannot
be a popping up isolated real node, nor two nodes can appear on the arc $a_5$.
Thus a possible position of the
new non-isolated node is determined by the rules (i) and (ii),
as we have seen in the base induction step.
\qed\enddemo

\proclaim{Remark \remone}{\rm
Statements similar to that of Lemma \lemone\ are contained also in
[\refsix, Step (1) in the proof of Lemma 3.3] and [\refthree, Lemma 2.10],
where detailed proofs have been skipped. Moreover, Lemma \lemone\ and the above cited
statements follow from [\refseven, Theorem 10 (proof in \S7-11)].
For the reader's convenience we have provided here a proof with all necessary
details that also complete the proofs in [\refsix, Step (1) in the proof of Lemma 3.3]
and [\refthree, Lemma 2.10].}
\endproclaim


\proclaim{Lemma \lemtwo}
Let $C_9$ be a real nodal sextic as in Lemma~\lemone and $p,q\in C_9$ be
as in Figure \figtwod, right.
Then there exists a real elliptic sextic $C(A_3)$ with
$7$ nodes and a singularity $A_3$ located with respect to the line $L$
as shown in Figure~\figtwob.
\endproclaim

\demo{Proof} We apply the Hilbert-Rohn-Gudkov method in the form developed in
[\refone, Section 4] and proceed similarly to the lines of the proof of
[\refone, Lemma 5.3].

\smallskip

By [\refthree, Proposition 2.4(b)], we can suppose that
$$
    \matrix
    \text{the configuration consisting of any prescribed $7$ nodes of $C_9$, }\\
    \text{the points $p$ and $q$, and of the tangent at $p$ is in general position.}
    \endmatrix
                                                                 \eqno(\eqtwo)
$$

Let us order the non-isolated nodes $z_1,\dots,z_4$ of $C_9$ assuming that $z_4\in\On$,
and respectively denote by $d'_1,\dots,d'_4$ the disks inside $d_0$ bounded by the arcs
ending at $z_1,\dots,z_4$ (so, $z_3,z_4\in d'_4$).
Pick a point $p\in a_6$ and a point $q\in\Oe$.
Consider the germ ${\Cal M}$ at $C_9$ of the equisingular family
of real elliptic sextics which
\roster
\item"(i)"
           have nodal singularities at all the isolated nodes of $C_9$ and
           at $z_1$ and $z_2$,
\item"(ii)"
           have a node in a neighborhood of $z_i$, $i=3,4$,
\item"(iii)"
           intersect $C_9$ at $p$ with multiplicity $2$,
\item"(iv)"
           pass through the point $q$.
\endroster

\midinsert
\epsfxsize=50mm
\centerline{
 \epsfbox{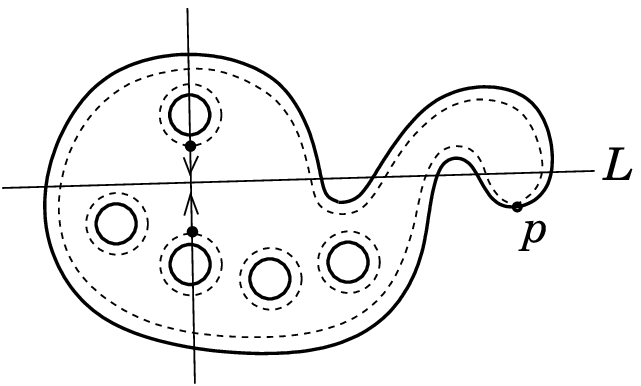}\hskip10mm
 \epsfxsize=50mm
 \epsfbox{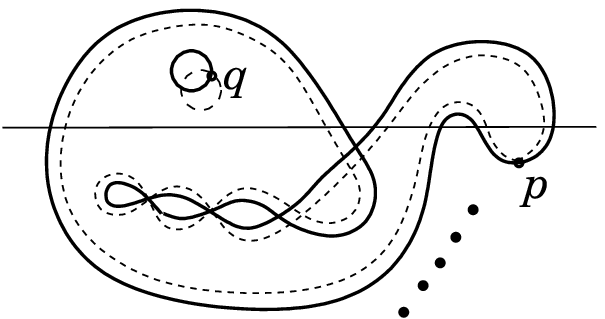}
}
\botcaption{ Figure \figtwod }
\endcaption
\endinsert

The germ ${\Cal M}$ is smooth and one-dimensional by [\reffour, Theorem in page 31].

By formula (24) in [\refone, Lemma 4.2] and formulas (15), (16) in
[\refone, Lemma 4.1], each curve $C'\in{\Cal M}\setminus\{C_9\}$ intersects $C_9$
with multiplicity $4$ at each of the seven fixed nodes, at two real points
in a neighborhood of $z_i$, $i=3,4$, and with multiplicity $3$ at $\{p,q\}$.
In total this gives $35$, and by the parity argument, one more (real) intersection point
of $C'$ with $C_9$ lies on the oval $\Oe$. Altogether this yields that, moving along
${\Cal M}$ in a certain direction, we obtain a deformation of the real point set such
that (see the dashed lines in Figure \figtwod, right):
$$
     \text{the disks $d'_1,d'_2,d'_3$ grow, the disk $d'_4$ and the domain $d_0$
     shrink;}
                                                      \eqno(\eqone)
$$
cf. [\refthree, Proposition 2.5 and Figure 3].

Extending the germ ${\Cal M}$ to a global equisingular family subject to
conditions (i), (iii), (iv) above, we see that the element of $\Cal M$
moving in the designated
direction cannot return to $C_9$ due to the strongly monotone changes (\eqone),
and hence must undergo a degeneration. The argument used in the proof of
Lemma~\lemone\ shows that it is not a double cubic.
The general position condition
(\eqtwo) excludes all other splittings of the degenerate curve into three or more
components (counting multiplicities).
Let us show that no splitting into two distinct components is possible.
Indeed, it follows from (\eqone), that the fixed isolated nodes in the domain $d_3$ remain isolated
in the degeneration, and that no component of odd degree can split off.
Note also that, in case of a splitting, no isolated node can be an intersection point of
two components, since otherwise, the components must be complex conjugate,
which is impossible. All this leaves the only possibility of a splitting into
a conic and a quartic, but such a curve cannot have $5$ isolated nodes.

Thus, the appearance of an extra node can only be the contraction of the oval $\Oe$
to the point $q$, otherwise, one would encounter a forbidden reducible curve.
In this case, we simply ignore the degeneration and continue
the movement along our one-dimensional family: the oval $O^e$ pop up again,
and the rest of the real part of the current curve deforms as shown in Figure \figtwod (right).
Then at some moment we have to encounter another degeneration, and the only possibility
left is the shrinking to a point of the disk $d'_4$, thus, giving the required
elliptic curve $C(A_3)$.
Indeed, the genus formula combined with (\eqone) do not allow
any singularity of the form $A_{n}$, $n\ne3$.
\qed\enddemo


\head \sectResolv. Application of cubic resolvents. End of proof of Theorem 1
\endhead

We denote the standard real Hirzebruch surface of degree $n>0$
(the fiberwise compactification
of the line bundle $\Cal O(n)$ over $\Bbb P^1$) by $\hirz_n$.
Let $\R\hirz_n$ be the set of real points of $\hirz_n$. It is diffeomorphic to
a torus or a Klein bottle. In Figures \figA\ and \figB\ we represent $\R\hirz_n$
by a rectangle whose opposite sides are identified.
The horizontal sides represent the exceptional section $E$, $E^2=-n$,
and vertical lines (in particular, the vertical sides of the rectangle) represent
fibers of the projection $\hirz_n\to\Bbb P^1$. Let $F$ be one of the fibers.
The Picard group of $\hirz_n$ is generated by $E$ and $F$.
A generic section disjoint from $E$ belongs to the linear system $|E+nF|$.

When speaking of a {\it fiberwise arrangement} of a curve on $\R\hirz_n$,
we mean its arrangement up to isotopies which fix $E$ and send each fiber to
a fiber. If it is known that the curve belongs to $|dE+ndF|$,
its {\it almost fiberwise arrangement} is the arrangement up to isotopies
fixing $E$ and such that any fiber at any moment intersects the curve at
$\le d$ points counting the multiplicities. In particular,
the ovals of trigonal curves ($d=3$) cannot pass one over another during
such isotopies.

\demo{ Proof of Theorem 1}
Suppose that Figure \figTh\ is realizable.
Then, by Lemma \lemtwo, there exists
a singular sextic curve with an $A_3$ singularity arranged with respect to $L$
as in Figure \figtwob. It can be perturbed into a curve
$C'$ arranged in one of the two ways shown in Figure \figA\ (left)
with respect to $L$ and the two dashed lines (by rotating $L$ around the
common point of the arcs $a_3$ and $a_6$ we can achieve that $L$ passes through $A_3$).
The relative position of $C'$ with respect to the (dashed) line through $\Oe$ follows
from [\refFO, Lemma 2.2]. The position of the two nodes
with respect to the tangent line at $A_3$ follows
from Bezout theorem for an auxiliary conic tangent to $C'$ at $A_3$ and passing through
$\Oe$, one empty exterior oval, and one of the nodes.

\midinsert
\epsfxsize=123mm\epsfbox{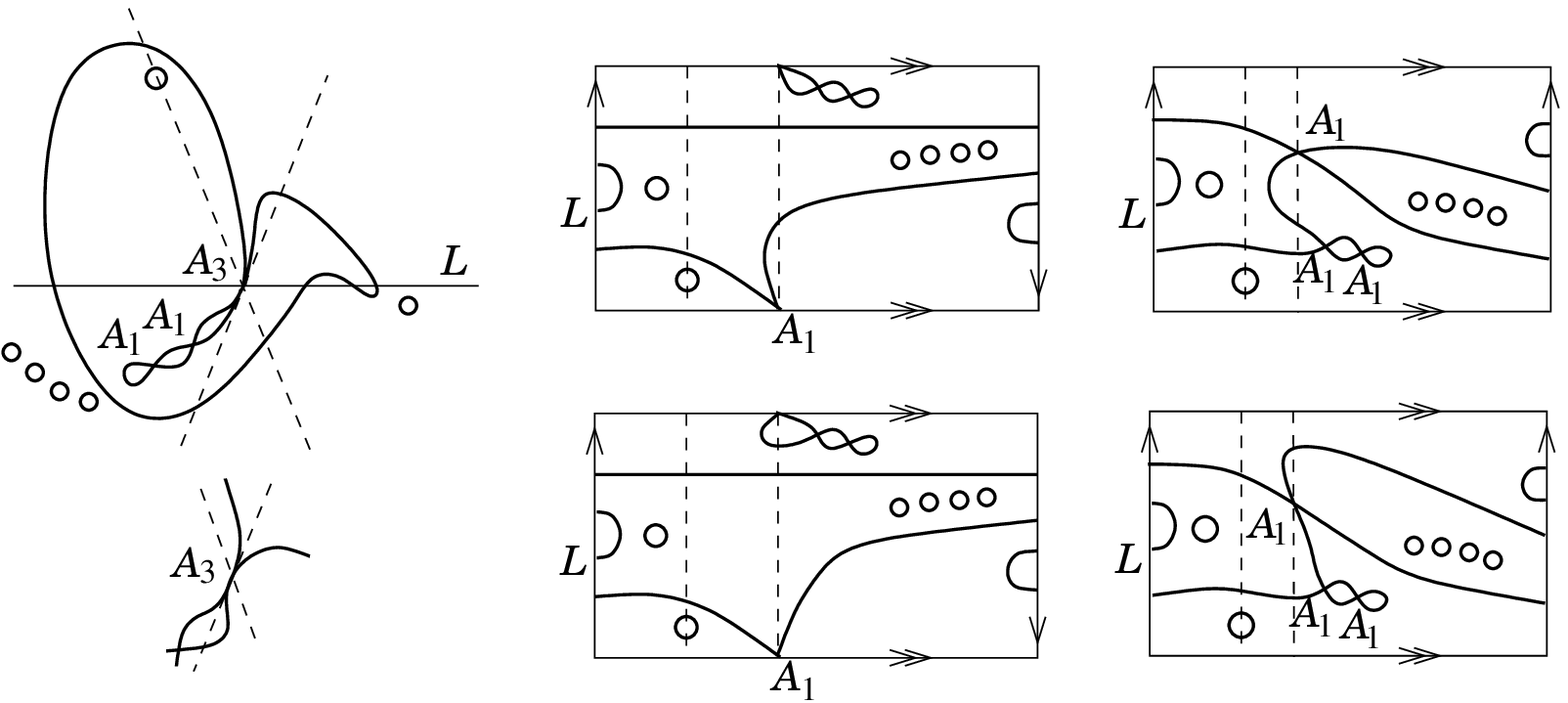}\vskip-3mm
\botcaption{ Figure \figA }
$C'\cup L$ and its transforms on $\Cal F_1$ and on $\Cal F_2$
\endcaption
\endinsert

Let us blow up the point $A_3$. We obtain the arrangement
in Figure \figA\ (middle) on $\R\hirz_1$.
Then we blow up the point $A_1$ on $E$ and blow down the strict transform of the
fiber passing through it. We obtain a curve in $\R\hirz_2$
belonging to $|4E+8F|$ arranged (up to isotopy) with respect to $E$ and the
indicated fibers as shown in Figure \figA\ (right).
Its cubic resolvent is a trigonal curve in $\hirz_4$ arranged
with respect to $E$ and the indicated fibers as in
Figure \figB\ (left); see [\refMMJ, \S3].

\midinsert
\epsfxsize=120mm\epsfbox{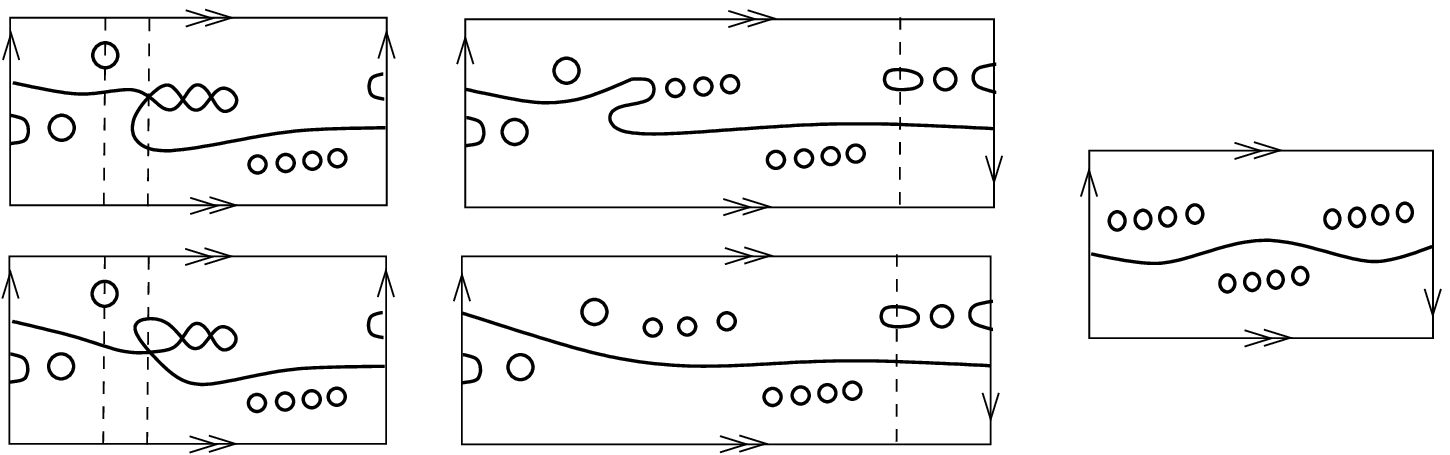}
\botcaption{ Figure \figB }
The cubic resolvent on $\Cal F_4$; the glued curve on $\Cal F_5$
\endcaption
\endinsert

Using [\refAFST], it is an easy exercise to check that the arrangement
in Figure \figB\ (left) is unrealizable by a trigonal algebraic curve on $\R\hirz_4$
(note that it is realizable by a trigonal pseudoholomorphic curve).
To this end one should exclude all its possible fiberwise arrangements,
namely, the one depicted in Figure \figB\ (left) and those obtained from it
by inserting a zigzag
\epsfysize=1.2ex\epsfbox{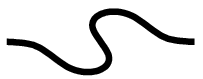} or
\epsfysize=1.2ex\epsfbox{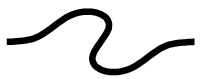}
between some two consecutive ovals
(at most one zigzag can be inserted because otherwise we
obtain too many vertical tangents).
Note that insertion of a zigzag is really necessary because there are
unrealizable fiberwise arrangements which become realizable after a
zigzag insertion, see [\refAa, Appendix B]; this phenomenon is impossible in
pseudoholomorphic context.

According to [\refAFST], to exclude each of these fiberwise arrangements, it
is enough to check that there does not exist a graph in $\CP^1$
satisfying Conditions (1)--(7) at the end
of [\refAFST, \S4] and having a prescribed behavior near $\RP^1$.
Indeed, since the number of vertices of the graph and their nature is
dictated by these conditions, only a finite number of cases should be considered
which can be done by hand in a reasonable time.
This fact can be also derived from Erwan Brugall\'e's result. He checked
in [\refB, Proof of Proposition 5.6] by this method that the almost fiberwise
arrangement in $\R\hirz_5$
shown in Figure \figB\ (right) is algebraically unrealizable. Indeed,
[\refB, Proposition 3.6] implies that it is enough to consider zigzag insertions
of the form \epsfysize=1.9ex\lower 0.3ex\hbox{\epsfbox{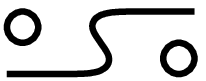}} up to symmetry, thus the
almost fiberwise unrealizability of Figure \figB\ (right) is a consequence
of [\refB, Lemmas 5.4 and 5.5].

The unrealizability of Figure \figB\ (left) follows from that of
Figure \figB\ (right) because the latter is obtained from the former
by gluing it together with an $M$-cubic in $\RP^2$ according to
Figure \figB\ (middle). The gluing can be understood either in the sense of [\refAFST]
or in the sense of Viro [\refV]. In the latter case we interpret the two
parts of Figure \figB\ (middle) as charts in the triangles
$[(0,0),(12,0),(0,3)]$ and $[(12,0),(15,0),(0,3)]$.
\qed\enddemo


\head \sectError. The mistake in [\refFO]
\endhead

The idea of the prohibition of the sextic in question realized in [\refFO]
 was to consider the pencil of real cubics through $8$ specific fixed  points on the hypothetical sextic, then, using an information on the  location of the fixed points with respect to lines and conic, to  construct the evolution of cubics along the pencil, and then to show  that such a pencil does not satisfy some necessary conditions (does  not reveal $8$ distinguished cubics, see definition in [\refFO, \S4]).

In this section we assume that the reader is familiar with the paper [\refFO] and we
use the notation from there.
The mistake is in the last step (``From $C^4$ to contradiction")
in [\refFO, \S5.4]: the assertion that the arcs 6 and 8 of $C^t$ are separated
by $C^4\cup N$ is erroneous.
As a matter of fact,
the non-base point of $C^t\cap N$ may escape
the loop of $N$ as shown in Figure~\figError.

\midinsert
\centerline{\epsfxsize=120mm\epsfbox{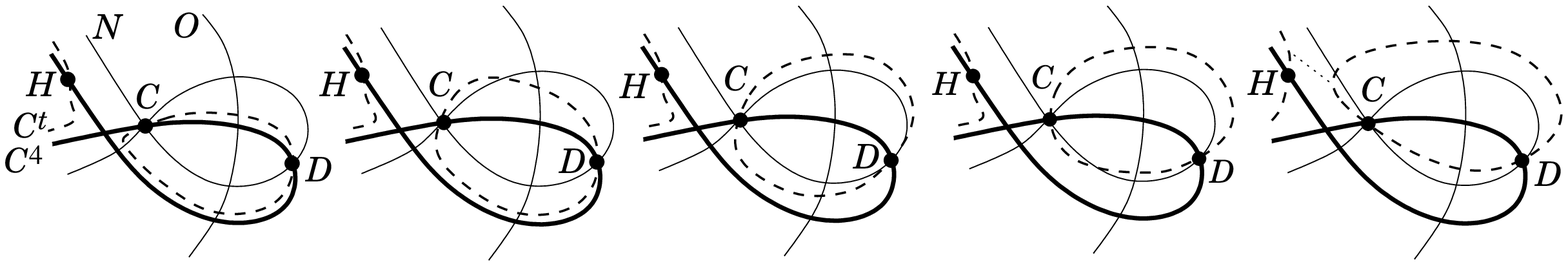}}
\botcaption{ Figure \figError } The non-base point of $C^t\cap N$ escapes the loop of $N$
\endcaption
\endinsert

One could hope to repair the proof in [\refFO] by continuing
the construction of the pencil of cubics and obtaining a contradiction on
some further step. Unfortunately, this is not so.
In Figure~\figPencil\ 
(see also Figure \figSev)
we complete the pencil of cubics without any
contradiction to Bezout theorem for the auxiliary curves considered in [\refFO].

\midinsert
\centerline{\epsfxsize=110mm\epsfbox{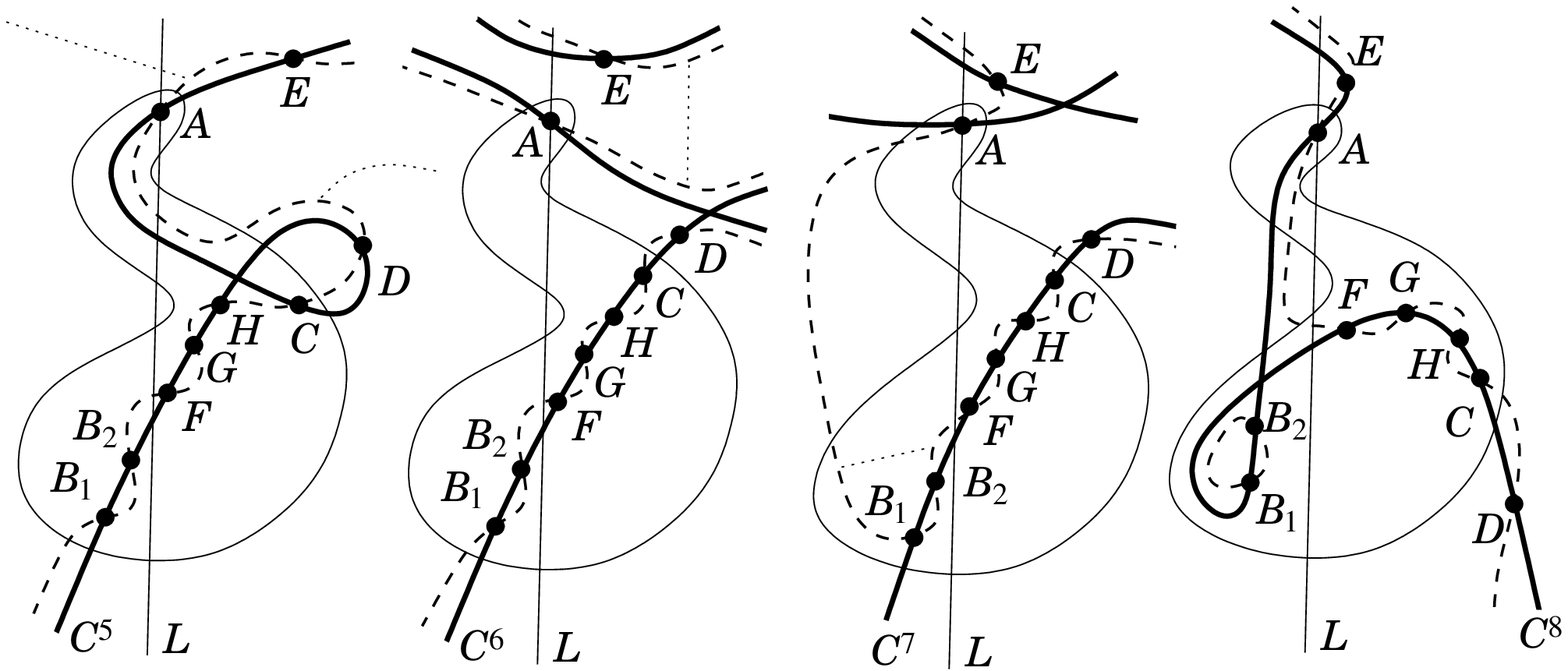}}
\botcaption{ Figure \figPencil } Completing the pencil of cubics (cf. [\refFO, Figure 17])
\endcaption
\endinsert


\if01{ 
The combinatorial structure of the considered pencil of cubics
(i.~e., the order of the base points along the distinguished cubics;
see [\refFO, \S4]) is uniquely determined by our setting. For the
first four distinguished cubics, this fact was established in [\refFO] by
a tedious case-by-case consideration. However the combinatorial uniqueness
of the whole pencil easily follows from the general results
[\refSevOne, \refSevTwo]. Indeed, ................
}\fi 

In the rest of the section we explain
how we construct the pencil, using the tools from [\refSevOne, \refSevTwo].
The {\it combinatorial configuration\/} of $n$ points in the plane is the data
describing the mutual position of each point with respect to the lines through
two others and the conics through five others. The {\em combinatorial pencil of cubics\/}
determined by eight points is given by the arrangement of the nine base points
on the eight successive {\it distinguished cubics} (see the definition in [\refFO, \S4]).
Let us consider eight points $(1, 2, 3, 4, 5, 6, 7, 8)$ distributed in the ovals
$(A, D, C, H, G, F, B, E)$. Using [\refFO, Lemma 2.2] plus Bezout's theorem
between $C_6$ and some auxiliary rational cubics (passing through seven of the points,
with node at one of them), we determine the combinatorial configuration $\Cal C$ realized by
these eight points. It is formed of five $7$-subconfigurations of type $(3, 4, 0, 0)_2$,
plus three of type $(7, 0, 0, 0)$ (see [\refSevTwo, \S3.1]).
To find this pencil determined by $1,\dots,8$, free the point $7$ away from the oval $B$ and move it till
it crosses the line $(AE)$. The new configuration $(1, \dots 8)$ lies in convex position,
its combinatorial configuration $\Cal C$ 
is replaced by $\max(\hat 1 = 8+)$ (see [\refSevOne, \S2.3]).
As $7$ is close to the line $(18)$, it lies outside of the loops of the cubics
$(\hat 7, 1)$ and $(\hat 7, 8)$, hence (see [\refSevOne, \S5.1])
the pencil is the first one in [\refSevOne, Figure 35].
Move $7$ back to its initial position in $B$,
the combinatorial pencil changes when $7$ crosses the line $(AE)$,
see upper part of Figure \figSev.
Afterwards, the eight points realize $\Cal C$ for all positions of $7$ on
the path. The only way to change the pencil would be to let $9$ cross another
base point $k$: when $9 = k$, the points $1,\dots,8$ lie on a nodal cubic with node at $k$
(see [\refSevOne]). But one proves that $\Cal C$ is not realizable by eight points on such a cubic.
So, the pencil undergoes no further change.
The complete pencil of cubics is shown in the lower part of this Figure \figSev.
(Note that with the notation of [\refFO], $B_1 = 9$ and $B_2 = 7$.)

\midinsert
\noindent\epsfxsize=130mm\epsfbox{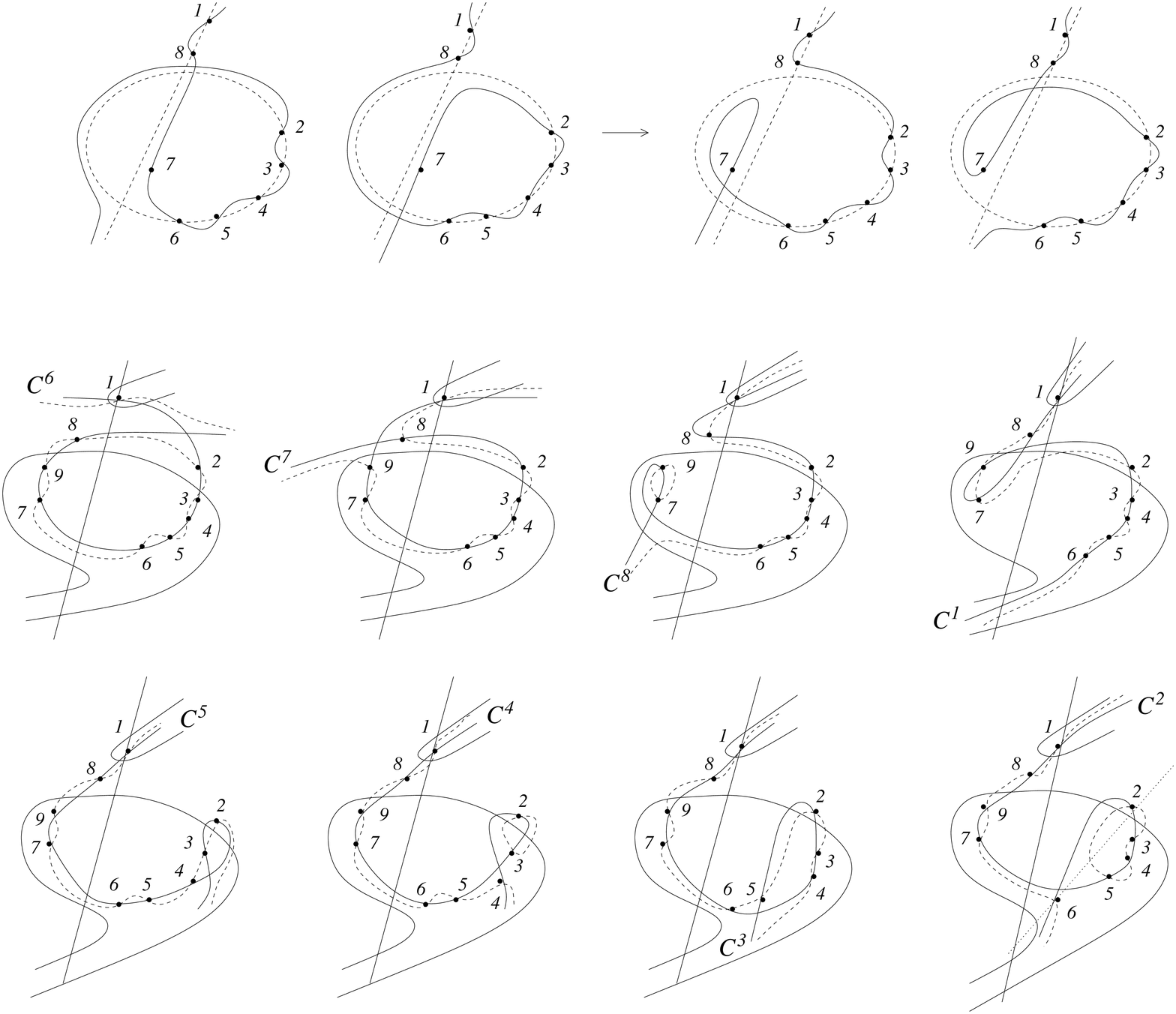}
\botcaption{Figure \figSev}
Upper part: change of two distinguished cubics induced letting $7$ cross the line $(18)$;
lower part: the complete pencil of cubics through points $1, 2, 3, 4, 5, 6, 7, 8$
distributed in the ovals $A, D, C, H, G, F, B, E$
\endcaption
\endinsert

\smallskip

{\bf Acknowledgements.} We are grateful to the referees for valuable comments and corrections.

\enddocument